\documentclass[12pt,preprint]{elsarticle}
\usepackage{graphicx}
\usepackage[english]{babel}
\usepackage{amsmath}
\usepackage{amssymb}
\usepackage{amstext}
\usepackage{lineno,hyperref}
\usepackage{caption}
\newdefinition{rmk}{Remark}
\begin{document}

\begin{frontmatter}

\title{Analytic Investigation for Spatio-temporal Patterns Propagation in Spiking Neural Networks}

\author{Ning Hua\corref{cor1}}
\ead{19110180033@fudan.edu.cn}
\address{School of Mathematical Sciences, Fudan University, Shanghai 200433, China}

\author{Xiangnan He}
\address{Centre for Computational Systems Biology, Fudan University, Shanghai 200433, China\\School of Mathematical Sciences, Fudan University, Shanghai 200433, China}

\author{Wenlian Lu\corref{cor2}}
\ead{wenlian@fudan.edu.cn}
\address{School of Mathematical Sciences, Fudan University, Shanghai 200433, China\\Shanghai Center for Mathematical Sciences, Fudan University, Shanghai 200438, China\\Shanghai Key Laboratory for Contemporary Applied Mathematics, Shanghai 200433, China}

\author{Jianfeng Feng}
\address{Institute of Science and Technology for Brain-Inspired  Intelligence, Fudan University, Shanghai 200433, China\\Key Laboratory of Computational Neuroscience and Brain-Inspired Intelligence (Fudan University), Ministry of Education, Shanghai 200433, China\\MOE Frontiers Center for Brain Science, Fudan University, Shanghai 200032, China}

\cortext[cor1]{The first author}

\cortext[cor2]{Corresponding author}

\begin{abstract}
Based upon the moment closure approach, a Gaussian random field is constructed to quantitatively and analytically characterize the dynamics of a random point field. The approach provides us with a theoretical tool to investigate synchronized spike propagation in a feedforward or recurrent spiking neural network. We show that the balance between the excitation and inhibition postsynaptic potentials is required for the occurrence of synfire chains. In particular, with a balanced network, the critical packet size of invasion and annihilation is observed. We also derive a sufficient analytic condition for the synchronization propagation in an asynchronous environment, which further allows us to disclose the possibility of spatial synaptic structure to sustain a stable synfire chain. Our findings are in good agreement with simulations and help us understand the propagation of spatio-temporal patterns in a random point field.
\end{abstract}


\end{frontmatter}

%\linenumbers
\section{Introduction}
Towards fully understanding an evolutionary random point field, the (joint) probability distribution density is often very hard, if it is not impossible, to be calculated and estimated. However, in many cases, the first few moments are sufficient to present a holistic picture of the evolution of this random point field\cite{Barzel12,Gleeson12,Adomian83,Xiao18}. To this end, the moment closure approach was proposed to approximate random fields by employing their first few moments \cite{MC1,MC2,MC3,MC4}, which transforms the spatio-temporal random point field into a few non-random dynamical systems of the moments respectively. Then, this idea enables us to investigate these moment dynamical systems for the depiction of the asymptotical properties of the random point field. For instance, the spatio-temporal pattern of random point field can be regarded as attracting dynamics of the corresponding moment dynamical system.

\begin{figure*}[!htb]
\centering\includegraphics[width=0.8\textwidth]{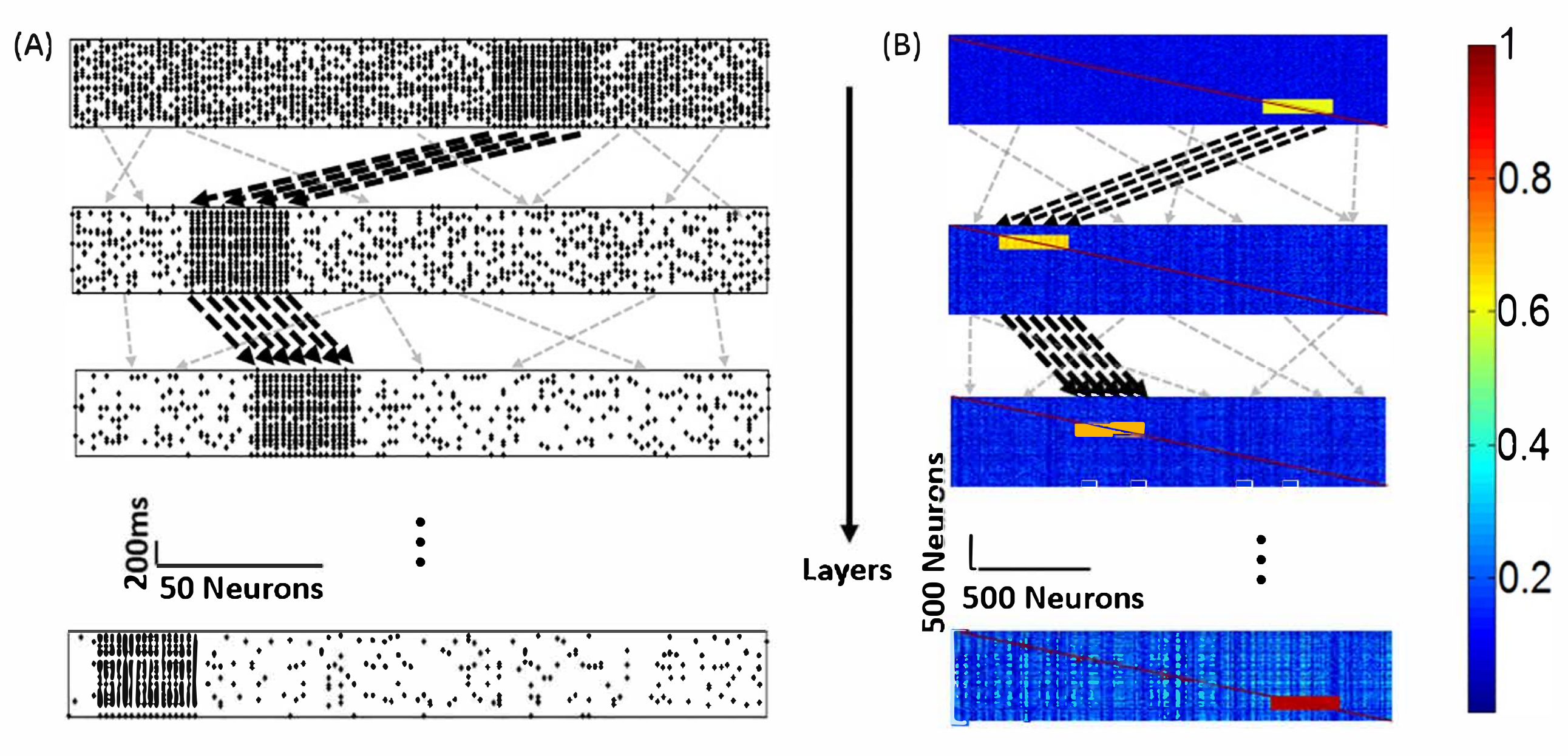}
\caption{Schematic plots of synfire chain. Spikes are transmitted from layer to layer (from top to bottom). Inside the pulse packet, each neuron is fully connected with neurons in the successive layer (dark dashed arrows); outside the pulse packet, the feedforward links are sparse (gray arrows).  (A). raster plots of spiking patterns (shadow spikes) in a feedforward spiking LIF neural network; (B). correlation matrices of synfire chain (warm-coloured rectangles) through layers. }
\label{fig1}
\end{figure*}

Mathematically, neural activities can be formulated as random point processes (spike trains), and a large ensemble of neural spiking trains with an underlying geometric structure
(manifold) composes of a random point field\cite{Wilson73}. In particular, the propagation of spatio-temporal patterns of neural spiking activities, as the essence of cortical function, has been widely studied in recent decades \cite{Mount,Kumar08,Cui07,Yanchuk11,Chiappalone05}.
A typical example of spatio-temporal pattern, the so-called pulse-packet, to encode a piece of information and reliably transmit it from one layer of nervous system to the other is
a synfire chain: at each layer, spikes are synchronized inside a packet of neurons but asynchronous out of the packet \cite{Abeles82,Riehle97,Jahnke2008,Aviel03} (see Fig. \ref{fig1} for illustration). There are large numbers of studies showing that synfire chain can transmit information between neuronal populations experimentally and theoretically\cite{Wang16,Shao16,Xiao17}. As an information transmitter, it has been successfully applied in many computational tasks in recent decades\cite{Jacquemin94,Abeles04,Hayon05,Izhikevich06}. This sort of coexistence of synchronous and asynchronous states of neural activities can also be categorised as "chimera" phenomenon \cite{chimera1,chimera2,chimera3}.

The synfire chain can be naturally described in the framework of random point fields. Traditionally, dynamical models of mean firing rate \cite{Shadlen94,Aertsen96} were used to depict the evolution of spike spatio-temporal patterns \cite{Adrian26,Barlow72,Wilson93,Romo99,Shadlen98} in cortical circuits. However, such models are only true under the condition of independent or weak correlation of spiking activities between neurons \cite{Mazurek02,Cateau01} and they have limitations to fully account for the evolution of the spiking patterns \cite{Mainen95,Prut98,Nowak97,Reyes03}. On the other hand, temporal dispersion in terms of  pulse-packets \cite{Dies} and correlation map \cite{Sompolinsky01,Rocha07} represented the degree of synchrony well but failed to realistically describe the pattern in a random point field. How the survival of a pulse packet in a multi-layered neural network depends on the structural and physiological characteristics of the network has not been answered analytically in the literature.

Under the framework of moment closure, Pearson correlation coefficients were utilized to describe synchrony \cite{Rocha07,Feng2006}. In the present paper, we developed a novel and general theoretical framework to investigate the synchrony dynamics of random point field of multilayer feedforward neural network (FNN) with spiking leaky integrate-and-fire neurons. Besides, we also find similar dynamical features in discrete-time recurrent spiking neural network using our model. By constructing the moment dynamical system that includes the first (mean firing rate) and second-order statistics (variance and correlation) of spiking random point fields, the synfire chain dynamics can be regarded as the existence of certain sort of attracting set of the correlation map. Combining with the mean field approach, we discover the necessity of the balance of network for the stable synchronization propagation. In addition, we analytically derive a sufficient condition for the existence of this attracting set that enables to obtain the appropriate size of this synfire chain and the proper synpatic density of the neural network. These results are in a good agreement with the numerical results and help understand the evolution of synchronisation pattern in spiking neural network.

\section{Methods}

\subsection{Gaussian approximation of random point field}
We start with an evolutionary random point field, $\xi(x,\tau)$ standing for the number of events occurring at $x$ for the duration $[0,\tau]$, where $x$ is the spatial variable and $t$ is time. We discretize the continuous-time into time bins $[0,\Delta,2\Delta,\cdots,k\Delta,\cdots]$ with a length $\Delta$. In each time bin, $[k\Delta, (k+1)\Delta]$ with $k\ge 0$, we define a random field $\xi^{k}(x,t)=\xi(x,t+k\Delta)-\xi(x,k\Delta)$, which counts the number of events occurring in $[k\Delta,k\Delta+t)$ for each $t\in[0,\Delta]$.

For each $\xi^{k}(x,t)$, the number of events occurring in an infinitesimal interval of time, $[t,t+dt]$, denoted by $d\xi^{k}(x,t)$, follows a probability distribution, $P^{k}(\xi;x,dt)$. Inspired by the moment closure approach, with assuming that $\xi(x,t)$ is a renewal process, which disregards the correlation with time, we are in the stage to approximate $P^{k}(\xi;x,dt)$ by the first moment, mean, the second statistics, variation, and correlation. Thus, we approximate $d\xi^{k}(x,t)$ by a Gaussian field as follows:
\begin{eqnarray}
d\xi^{k}(x,t)\approx\mu^{k}(x,t)dt+\sigma^{k}(x,t)d\eta^{k}(x,t)\label{variation}
\end{eqnarray}
where $\mu^{k}(x,t)$ and $\sigma^{k}(x,t)$ are the functions of mean and variance, $\eta^{k}(x,t)$ is a Gaussian field with zero means and unit variances: $E(\eta^{k}(x,t))=0$, $var(d\eta^{k}(x,t))=dt$, but correlated. Here, $E(\cdot)$ stands for the expectation. Equivalently, we can use the following coefficient of variation (CV) to stand for the second moment:
\begin{eqnarray}
CV^{k}(x,t)=\frac{\sigma^{k}(x,t)}{\mu^{k}(x,t)}.\label{CV}
\end{eqnarray}
The Pearson correlation coefficient (CC) between locations $x$ and $y$ at $t$:
\begin{eqnarray}
\rho^{k}(x,y,t)=\frac{E\left[d\eta^{k}(x,t),d\eta^{k}(y,t)\right]}{dt}.\label{CC1}
\end{eqnarray}

Assume that $d\xi^{k}(x,t)$ is asymptotically stationary. Then, with a sufficiently large $\Delta$, the long-term duration behaviours of $\mu^{k}(x,t)$, $\sigma^{k}(x,t)$ and $\rho^{k}(x,y,t)$ can be represented by the asymptotics, namely, $\mu^{k}(x)$, $\sigma^{k}(x)$ and $\rho^{k}(x,y)$. Thus, (\ref{variation}) can be rewritten as:
\begin{eqnarray}
d\xi^{k}(x,t)\approx\mu^{k}(x)dt+\sigma^{k}(x)d\eta^{k}(x,t).\label{variation1}
\end{eqnarray}
Here, the correlation between $d\eta^{k}(x,t)$ and $d\eta^{k}(y,t)$ is $\rho^{k}(x,y)dt$. In this stationary case, CC can be equivalently defined in the following way. The shift-correlation between $\xi^{k}(x,t)$ and $\xi^{k}(y,t)$ in a sliding window with length $T$ is defined as:
\begin{eqnarray}\label{rho_T}
\rho_{T}^{k}(x,y)=\frac{cov(n^{k}(x,\tau),n^{k}(y,\tau))}{\sqrt{var(n^{k}(x,\tau))var(n^{k}(y,\tau))}},
\end{eqnarray}
where $n^{k}(x,\tau)$ is the number of events occurring at location $x$ in the time interval $[\tau,\tau+T]$, $cov(\cdot,\cdot)$ stands for the covariance and the covariance ($cov$) and variance ($var$) are calculated with respect to $\tau$. Thus, the correlation coefficient between this spike train pair is defined as its
limit $\rho^{k}(x,y)=\lim_{T\to\infty}\rho_{T}^{k}(x,y)$, which equals formula (\ref{CC1}) in stationary state\cite{Rocha07}.
\subsection{Moment clourse method}
The above analysis constructs a random point field system (\ref{variation1})-(\ref{rho_T}) to study the dynamical characteristics of spiking neural network. Discribing spike activities with probability distribution function directly seems intractable, so it is more efficient to focus on their moments of distribution. Specifically, suppose $\Lambda$ is a finite index set and $\mathbb{M}\triangleq\big\{m_{i}:i\in\Lambda\big\}$ represents the moments we need. Then, we can transfer a stochastic system into moment system, $\big\{m_{i}:i\in\Lambda\big\}=F(\big\{m_{i}:i\in\Lambda\big\})$, using moment closure methods \cite{MC4}, where $F(\cdot)$ represents corresponding moment map. Different from the famous model proposed by Wilson and Cowan \cite{wca} which only includes the first order moment (firing rate), we hope to find a more complex model based on higher order, such as the second order moments (variance, correlation). Gaussian part in the right-hand side of (\ref{variation1}) can naturally provide the informatin of the second order moments and yield a moment closure.
 
We regard the long-term dynamics in the time bin $[(k-1)\Delta, k\Delta]$, $\xi^{k-1}(x,t)$, as the input to that of the successive time bin, $[k\Delta, (k+1)\Delta]$,  $\xi^{k}(x,t)$. We aim to formulate the iteration of these first and second moments/statistic in the form of
\begin{eqnarray}
\begin{cases}
\mu^{k}(x)=\mathcal S_{1}^{x}\left[\mu^{k-1}(\cdot),cov^{k-1}(\cdot,\cdot)\right],\\
\sigma^{k}(x)=\mathcal S_{2}^{x}\left[\mu^{k-1}(\cdot),cov^{k-1}(\cdot,\cdot)\right],\\
\rho^{k}(x,y)=\Psi^{x,y}\left[\mu^{k-1}(\cdot),cov^{k-1}(\cdot,\cdot)\right],
\end{cases}\label{MNN1}
\end{eqnarray}
where $cov^{k-1}(x,y)=\rho^{k-1}(x,y)\sigma^{k-1}(x)\sigma^{k-1}(y)$ is the covariance of $\xi^{k-1}(\cdot,t)$ between the location $x$ and $y$. Here, $\mathcal S^{x}_{1,2}$ stand for the functions of mean and variance with respect to the last moment mean and covariance, named mean and variance maps respectively. $\Psi^{x,y}$ stands for the correlation function, named correlation map. The details of $\mathcal S^{x}_{1,2}$ and $\Psi^{x,y}$ will be introduced in Section 2.3. Thus, the set of  mean, variance and correlation maps compose of the moment dynamic system (\ref{MNN1}) that represents the random point field $\xi(x,\tau)$ in the moment closure fashion.

The problem of spatio-temporal pattern of $\xi(x,\tau)$, which can be described by the first and second moments/statistics, can be transformed into the dynamics of (\ref{MNN1}). For instance, synchronization pattern of a random point process can be described by the CC (\ref{CC1}). $\rho^{k}(x,y)=1$ implies the processes at $x$ and $y$ completely synchronise; $\rho^{k}(x,y)=0$ means uncorrelated and $\rho^{k}(x,y)=-1$ means complete desynchronization. The asymptotic dynamics, namely, asymptotic attractor, can depict the pattern related to synchronization. If $\xi(x,t)$ is asymptotically steady, we can equivalently study the equilibrium functions: $\mu^{k-1}(x)=\mu^{k}(x)$, $\sigma^{k-1}(x)=\sigma^{k}(x)$ and $\rho^{k-1}(x,y)=\rho^{k}(x,y)$, towards understanding the asymptotic pattern of the random point process. We highlight that this approach can be extended to include the high-order moments and correlations of the random point field.

\subsection{Multilayer feedforward spiking network}
As an application, let us specify a random point field described by multilayer FNN of spiking neurons with a sparse random coupling structure \cite{Aviel03,Kumar10}. At each layer, there are exactly $N=5n$ neurons. Among them, there are $N_E=4n$ excitatory (E-) neurons and $N_I=n$ inhibitory (I-) neurons. For the E-neuron group at the $k$th layer, a synfire chain is carried by a packet of E-neurons, denoted by $\mathcal W^{k}$ with the identical size $\# \{\mathcal W^{k}\}=W$ for all $k$, as the packet size, where $\# \{\cdot\}$ is the number of elements in a finite set. Between any two successive layers, for instance, from the $(k-1)$th to the $k$th layer, each neuron in $\mathcal W^{k}$ receives inputs from all neurons in $\mathcal W^{k-1}$. This constructs a completely connected FNN in the packet, which is believed as the basic model for synfire chain \cite{Kumar10}. This completely connected FNN is embedded into a sparse network \cite{Aviel03} and the couplings for other neurons outside $\mathcal W^{k}$ are randomly picked with an equal probability so that the total number of excitation and inhibition synaptic links are $K_E=\lambda N_{E}$ and $K_I=\lambda N_{I}$ respectively (see Fig. \ref{fig1}). Here, $\lambda\in (0,1)$ is the indicator of the  sparsity of the synaptic density. To ensure that every neuron receives an equal number of synapses, $W \leqslant K_{E}$ is necessary.

We use a discretized random point field $N^{k}_{i}(t)$ to stand for the spike counts from the neuron $i$ at the $k$th layer, i.e., $N^{k}_{i}(t)=\sum_{l}\delta(t-t^{i,k}_{l})$, where $\delta(\cdot)$ is the Kronecker-delta function and $\{t^{i,k}_{l}\}$ are the time points of the pre-synaptic spikes from neuron $i$ (at the $k$th layer). Here, the neuron label $i$ stands for the spatial variable, $t$ for the continuous-time variable at $k$th layer.

The evolution of the random point field follows the integration-and-fire model. The potential activity of each neuron $i$ at layer $k$ is described as:
\begin{eqnarray}\label{LIF}
\tau_{m}dV_{i}^{k}(t)=- V_{i}^{k}(t) dt + I_{ext,i}^{k}+I_{syn,i}^{k},
\end{eqnarray}
for $i=1,\cdots,N$. Here, $\tau_{m}$ is the capacitance constant, $I_{ext,i}^{k}$ is the external current stimulus at layer $k$, and
$I_{syn,i}^{k}$ is the synaptic stimulus from the neurons at the $(k-1)$th layer:
\begin{eqnarray}
I_{syn,i}^{k}=\sum_{j}w_{ij}^{k}dN_{j}^{k-1}\label{syn_stm}
\end{eqnarray}
with $w^{k}_{ij}$ standing for the strength of the excitatory post-synaptic potential (EPSP) or inhibitory post-synaptic potential (IPSP) from neuron j at the $(k-1)$th layer to neuron $i$ at the $k$th layer. In this paper, we take values of EPSP/IPSP as follows. If there is a synaptic link from neuron $j$ at the $(k-1)$th layer to neuron $i$ at the $k$th layer, then we set $w_{ij}^{k}=w_{0}$ if $j$ is an E-neuron and $w_{ij}^{k}=-rgw_{0}$ if $j$ is an $I$-neuron, for some constant $w_{0}>0$, where $g=4$ equalizes the ratio between the numbers of E-synapses over I-synapses and thus $r$ serves as the ratio of ISPS over ESPS. If there is no link from $j$ to $i$, $w_{ij}=0$.

Once $V_{i}^{k}(t)$ reaches a threshold ($V_{th}$), neuron $i$ is depolarized and emits a spike. Then, $V_{i}^{k}(t)$ is reset to $V_{r}$ after a period of refractory time $T_{ref}$. By this way, neuron $i$ at layer $k$ emits a spike train, which is the input to the neurons at the $(k+1)$th layer that are linked with neuron $i$. Therefore, the random point field of spiking trains, $N^{k}_{i}(t)$, can be generated.

By the approach mentioned above, we are to approximate the spiking train of neuron $j$ (a random point process) by a Gaussian process: $d N^{k-1}_{j}(t) \sim \mu^{k-1}_jdt+\sqrt{\tau_{m}}\sigma^{k-1}_jdB^{k-1}_j$, where $B^{k-1}_{j}$, $j=1, \cdots, N$ are correlated Brownian motions \cite{Feng2006,Lu2010}. $\mu^{k-1}_j$ and $\sigma^{k-1}_j$ are the mean and variance of neuron $j$ (at the $(k-1)$th layer), which can be derived by the renew process theory \cite{Cox}. Here, $\tau_{m}$ stands for the time-scale constant. According to (\ref{syn_stm}), (\ref{LIF}) becomes correlated Ornstein-Uhlenbeck (OU) processes:
\begin{equation}\label{OUapproximation}
 \tau_{m}dV_{i}^{k}(t)=-V_{i}^{k}(t)dt+\hat{\mu}_{i}^{k}dt+\sqrt{\tau_{m}}\hat{\sigma}_{i}^{k}
 dB_{i}^{k}(t),
\end{equation}
$i=1,\cdots,N$. Here $\hat{\mu}_{i}^{k}=\sum_{j}w^{k}_{ij}\mu^{k-1}_{j}$ and $\hat{\sigma}_{i}^{k}=\sqrt{\sum_{j,l}w_{ij}^{k}\sigma_{j}^{k-1} \rho_{jl}^{k-1}w_{il}^{k}\sigma_{l}^{k-1}}$ are the mean and variance of the sum of the post-synapses received by neuron $i$ respectively, where $\rho_{jl}^{k-1}$ is the correlation coefficient between the $j$-th and $l$-th synapses of neuron $i$.

The output spike trains derived from (\ref{OUapproximation}) (with potential threshold $V_{th}$) are also approximated as Gaussian processes. This establishes maps of the first (mean map), second-order statistics (variance and correlation maps) of the random point field $N^{k}_{i}(t)$ of the successive layers. By the Siegert's expression \cite{Feng2003}, these three maps can be formulated:
\begin{eqnarray}\label{MNN}
\mu_{i}^{k}&=&\mathcal{S}_1(\hat{\mu}_{i}^{k},\hat{\sigma}_{i}^{k})\nonumber\\
\sigma_{i}^{k}&=&\mathcal{S}_2(\hat{\mu}_{i}^{k},\hat{\sigma}_{i}^{k})
\sqrt{\mathcal{S}_1(\hat{\mu}_{i}^{k},
\hat{\sigma}_{i}^{k})},i=1, \cdots, N,
\end{eqnarray}
where $\mathcal{S}_1$ and $\mathcal{S}_2$ are the mean and variance maps. Specifically, they can be written as\cite{Feng2006,Lu2010}:
\begin{equation}
\begin{aligned}
\mathcal{S}_{1}(y, z) & \approx \frac{1}{\left(T_{r e f}+\frac{2}{L} \int_{I\left(V_{r}, y, z\right)}^{I\left(V_{t h}, y, z\right)} D_{-}(u) d u\right)}, \\
\mathcal{S}_{2}(y, z) & \approx \frac{\left(\frac{8}{L^{2}} \int_{I\left(V_{r}, y, z\right)}^{I\left(V_{t h}, y, z\right)} D_{-} \otimes D_{-}(u) d u\right)^{1 / 2}}{\left(T_{r e f}+\frac{2}{L} \int_{I\left(V_{r}, y, z\right)}^{I\left(V_{t h}, y, z\right)} D_{-}(u) d u\right)}.
\end{aligned}\label{momentmap}
\end{equation}
Here $T_{ref}$ represents the refractory period, $V_{r}$ is the rest potential and $V_{th}$ is the threshold of membrane potential to emit a spike. Besides, 
\begin{equation}
\begin{aligned}
I(\xi, y, z) &=\frac{\xi L-y}{z}, \\
D_{-}(u) &=\exp \left(u^{2}\right) \int_{-\infty}^{u} \exp \left(-v^{2}\right) d v, \\
D_{-} \otimes D_{-}(u) &=\exp \left(u^{2}\right) \int_{-\infty}^{u} \exp \left(-v^{2}\right) D_{-}^{2}(v) d v,
\end{aligned}
\end{equation}
where $D_-(u)$ is exactly the Dawson's integral. For the details, we refer the readers to Ref.~\cite{Feng2006,Lu2010}.

As pointed out in Ref.~\cite{Rocha07}, the evolution of correlation coefficient from the $(k-1)$th layer to the $k$th is formulated as:
\begin{eqnarray}
\rho_{i,j}^{k}&=\Phi(\hat{\rho}_{i,j}^{k}),~i,j=1, \cdots, N,\label{cc}
\end{eqnarray}
where
$\hat{\rho}_{i,j}^{k}$ is the correlation coefficient between the collection of synpatic inputs of neuron $i$ and $j$ at layer $k$:
\begin{eqnarray}\label{Rho}
\hat{\rho}_{i,j}^{k}=\frac{\sum\limits_{p,q}w_{ip}^{k}\sigma_{p}^{k-1}w_{jq}^{k}
\sigma_{q}^{k-1}\rho^{k-1}_{pq}}
{\hat{\sigma}_{i}^{k}\hat{\sigma}_{j}^{k}},\nonumber\\
p,q=1,2,\dots,K_E+K_I.\qquad
\end{eqnarray}
$\Phi$ is the correlation map, which is assumed to be the identity map, namely, $\Phi(\rho)=\rho$ in this paper, according to the arguments in Ref.~\cite{Feng2006}. However, the following results still hold when $\Phi(\mu)$ is monotonous increasing with respect to $\mu$ as shown in Ref.~\cite{Rocha07,Lu2010} with some minor modifications.

To sum up, the spike trains of neurons in the FNN are modeled by the discrete random point field $N^{k}_{i}(t)$. Their evolution through layers is formulated by the iteration of the first and second moments/statistics, which is completely described by the mean, variance,  correlation maps and called the moment neural network (MNN)\cite{Feng2006}.

In particular, evolution equations of correlation map through a multilayer FNN can be utilized to analytically and quantitatively study the synchronization propagation. The (stable) synfire chain is defined by an attractor of the correlation map (\ref{cc}) together with the mean and variance maps (\ref{MNN}) with two properties: (${\mathbf P_{1}}$)
the correlation coefficients between the neurons in the packet of synchrony are large; (${\mathbf P_{2}}$) the  correlation coefficients between the other pairs of neurons are low.

Hence to characterize the synfire chain, we focus on the correlation map, which turns out to be a useful way to understand the complex dynamics.
Let us use the subscript $+$ to denote all variables in the packet $\mathcal W^{k}$ and $-$ out the packet at each layer.  Then, for given two small constants $0< \epsilon<1$ and $0<\delta<1$, define $
\Lambda_{\epsilon,\delta}=\{(\rho_+,\rho_-):~\rho_{+}>1-\epsilon,~|\rho_{-}|<\delta\}
$. If $\Lambda_{\epsilon,\delta}$ is an attracting set of the feedforward network of LIF neurons or the theoretical model (\ref{cc}), then we
call it {\it synfire attractor} and $\rho_+-\rho_-$ the {\it synfire gap}. Thus, a synfire chain can be defined as:
\begin{rmk}
The multilayer FNN is said to possess a synfire chain (with respect to $\epsilon$ and $\delta$) if $\Lambda_{\epsilon,\delta}$ is an attracting set of the iteration map (\ref{cc}) accompanied with (\ref{MNN}).
\end{rmk}
\subsection{Spatial mean-field approximation}
With the mean field approximation, the correlation map allows to provide analytic inference and insights on the stability of synfire chain by some algebras, and so greatly simplifies the analysis of the transmission of spiking point process in multilayer FNN. According to whether the link or neuron belongs to the packet or not, we substitute the specific correlation between neurons $p$ and $q$, $\rho_{pq}^{k}$, and the variance of neuron $p$, $\sigma_{p}^{k}$, by $\rho_{\pm}^{k}$ and $\sigma_{\pm}^{k}$. In detail, the mean in-packet can be calculated through (\ref{cc})-(\ref{Rho}) as:
\begin{eqnarray}\label{rho_in}
\rho_{+}^{k+1}&=&\left\langle\Phi(\hat{\rho}_{i,j}^{k+1})\right\rangle_{(i,j)\in\mathcal{W}^{k+1}\times\mathcal{W}^{k+1}}\nonumber\\
&\approx &\frac{A_{+}^{k+1}}{B_{+}^{k+1}}.
\end{eqnarray}
Here $\left\langle\cdot\right\rangle_S$ represents the average calculation over the set $S$, noting that the correlation is assumed to be the identity map in this paper.
Similarly, the out-of-packet correlation coefficients are approximately written as:
\begin{eqnarray}\label{rho_out}
\rho_{-}^{k+1}&=&\left\langle\Phi(\hat{\rho}_{i,j}^{k+1})\right\rangle_{(i,j)\not\in\mathcal{W}^{k+1}\times\mathcal{W}^{k+1}}\nonumber\\
&\approx &a_1\bigg(\frac{A_{-}^{k+1}}{B_{-}^{k+1}}\bigg)+a_2\bigg(\frac{A_{-}^{k+1}}
{\sqrt{B_{+}^{k+1}B_{-}^{k+1}}}\bigg)\approx\frac{A_{-}^{k+1}}{B_{-}^{k+1}},\nonumber\\
\end{eqnarray}
where $a_1=\frac{(N-W)^2}{N^2-W^2}\approx 1$ and $a_2=\frac{2W(N-W)}{N^2-W^2}\approx 0$, when $N>>W$, are two weights corresponding to the proportions of neuron pairs with all out-of-packet and neuron pairs with one in-packet and the other out-of-packet respectively. Here,
\begin{eqnarray}
\begin{cases}
A_{+}^{k+1}=\left\langle\sum_{p,q}w_{ip}^{k+1}\sigma_{p}^kw_{jq}^{k+1}\sigma_{q}^k\rho^{k}_{pq}\right\rangle_{i,j\in\mathcal W^{k+1}},\\
B_{+}^{k+1}=\left\langle\sum_{p,q}w_{ip}^{k+1}\sigma_{p}^kw_{iq}^{k+1}\sigma_{q}^k\rho^{k}_{pq}\right\rangle_{i\in\mathcal W^{k+1}},\\
A_{-}^{k+1}=\left\langle\sum_{p,q}w_{ip}^{k+1}\sigma_{p}^kw_{jq}^{k+1}\sigma_{q}^k\rho^{k}_{pq}\right\rangle_{i,j\notin\mathcal W^{k+1}},\\
B_{-}^{k+1}=\left\langle\sum_{p,q}w_{ip}^{k+1}\sigma_{p}^kw_{iq}^{k+1}\sigma_{q}^k\rho^{k}_{pq}\right\rangle_{i\notin\mathcal W^{k+1}}.\label{expofAB}
\end{cases}
\end{eqnarray}

Through a squared decomposition in Fig. \ref{fig2}, we have the following expression of $A_{\pm}^{k+1}$ and $B_{\pm}^{k+1}$ in mean field form. See Appendix A for the detailed derivation.
\begin{small}
\begin{align}\label{Variance_11}
&A_{+}^{k+1}\approx W(\sigma_{+}^k)^2(1-\rho_{+}^k)+\lambda(K_{E}-W)(\sigma_{-}^k)^2(1-\rho_{-}^k)+4r^2
\lambda{K_{E}}(\sigma_{-}^k)^2(1-\rho_{-}^k)+\nonumber\\
&\rho_{-}^k\!\bigg\{\!(r\!-\!1)^2\ K_{E}^2(\sigma_{-}^k)^2\!+\!2K_{E}W\sigma_{-}^k(\sigma_{+}^k\!-\!\sigma_{-}^k)
\!+\!W^2(\sigma_{+}^k\!-\!\sigma_{-}^k)^2\!\bigg\}\!\!+\!W^2(\sigma_{+}^k)^2(\rho_{+}^k\!-\!\rho_{-}^k),
\end{align}
\begin{align}\label{Variance_12}
&B_{+}^{k+1}\approx
W(\sigma_{+}^k)^2(1-\rho_{+}^k)+(K_{E}-W)(\sigma_{-}^k)^2(1-\rho_{-}^k)
+4r^2{K_{E}}(\sigma_{-}^k)^2(1-\rho_{-}^k)+\nonumber\\
&\rho_{-}^k\!\bigg\{\!(r-1)^2K_{E}^2(\sigma_{-}^k)^2
\!+\!2K_{E}W\sigma_{-}^k(\sigma_{+}^k\!-\!\sigma_{-}^k)\!+\!W^2(\sigma_{+}^k\!-\!\sigma_{-}^k)^2\!\bigg\}\!\!+\!W^2(\sigma_{+}^k)^2(\rho_{+}^k\!-\!\rho_{-}^k),
\end{align}
\begin{align}\label{Variance_13}
&A_{-}^{k+1}\approx\lambda^2W(\sigma_{+}^k)^2(1\!-\!\rho_{+}^k)+\lambda(K_{E}\!-\!\lambda{W})(\sigma_{-}^k)^2(1\!-\!\rho_{-}^k)
\!+\!4r^2\lambda{K}_{E}(\sigma_{-}^k)^2(1\!-\!\rho_{-}^k)\!+\!\nonumber\\
&\rho_{-}^k\bigg\{(r-1)^2K_{E}^2(\sigma_{-}^k)^2
\!+\!2\lambda{K}_{E}W\sigma_{-}^k(\sigma_{+}^k\!-\!\sigma_{-}^k)\!+\!\lambda^2W^2(\sigma_{+}^k-\sigma_{-}^k)^2\bigg\}+\lambda^2W^2(\sigma_{+}^k)^2\nonumber\\&(\rho_{+}^k-\rho_{-}^k),
\end{align}
\begin{align}\label{Variance_14}
&B_{-}^{k+1}\approx\lambda{W}(\sigma_{+}^k)^2(1-\rho_{+}^k)+(K_{E}-\lambda{W})(\sigma_{-}^k)^2(1-\rho_{-}^k)
+4r^2{K_{E}}(\sigma_{-}^k)^2(1-\rho_{-}^k)+\nonumber\\
&\rho_{-}^k\!\bigg\{\!(r\!-\!1)^2\lambda^2{N}_{E}^2(\sigma_{-}^k)^2
\!+\!2\lambda{K}_{E}W\sigma_{-}^k(\sigma_{+}^k\!-\!\sigma_{-}^k)
+\lambda^2{W}^2(\sigma_{+}^k-\sigma_{-}^k)^2\bigg\}+\lambda^2{W}^2(\sigma_{+}^k)^2\nonumber\\&(\rho_{+}^k-\rho_{-}^k).
\end{align}
\end{small} 
\begin{figure*}[!htb]
\centering\includegraphics[width=0.8\textwidth]{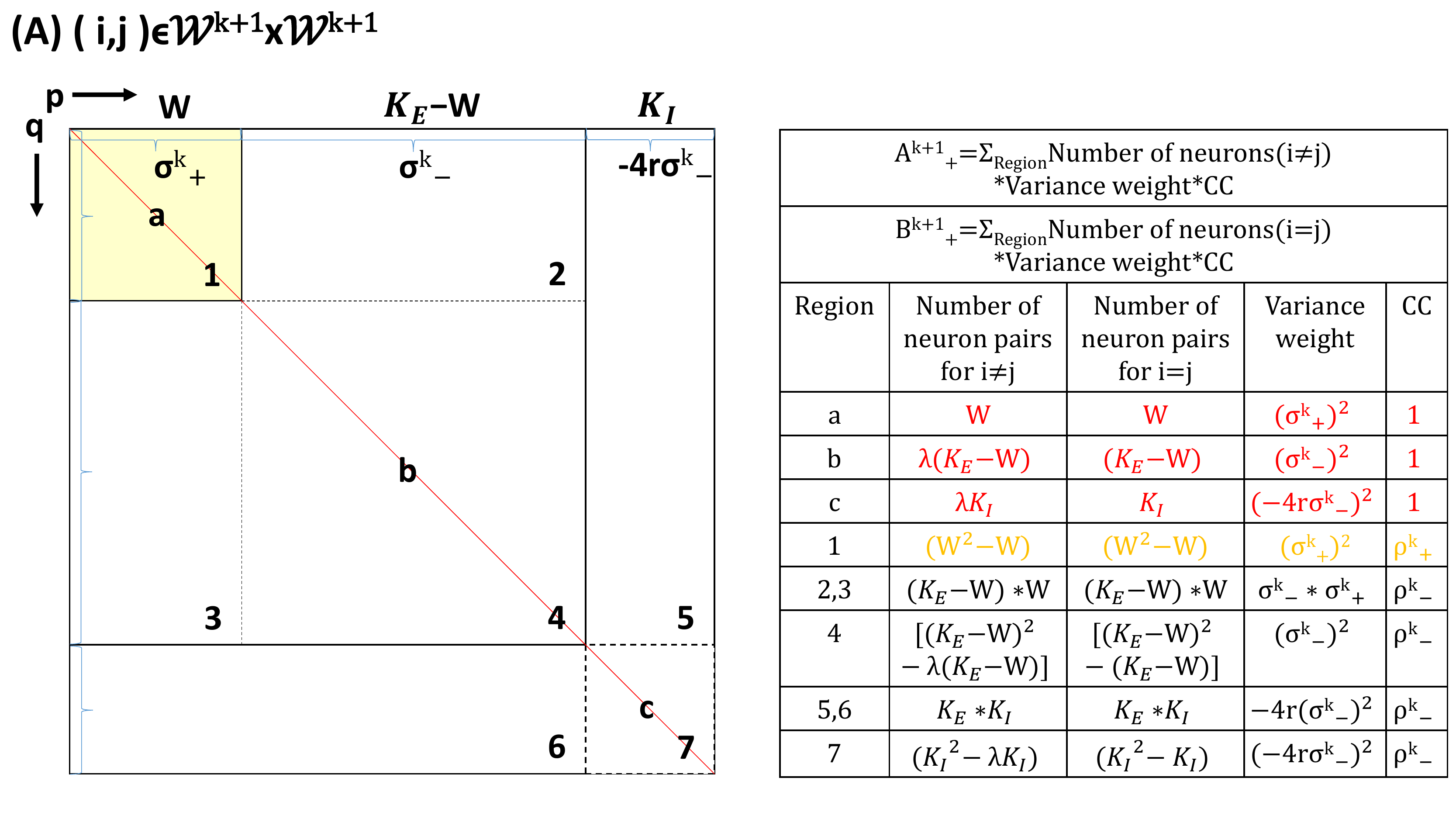}
\centering\includegraphics[width=0.8\textwidth]{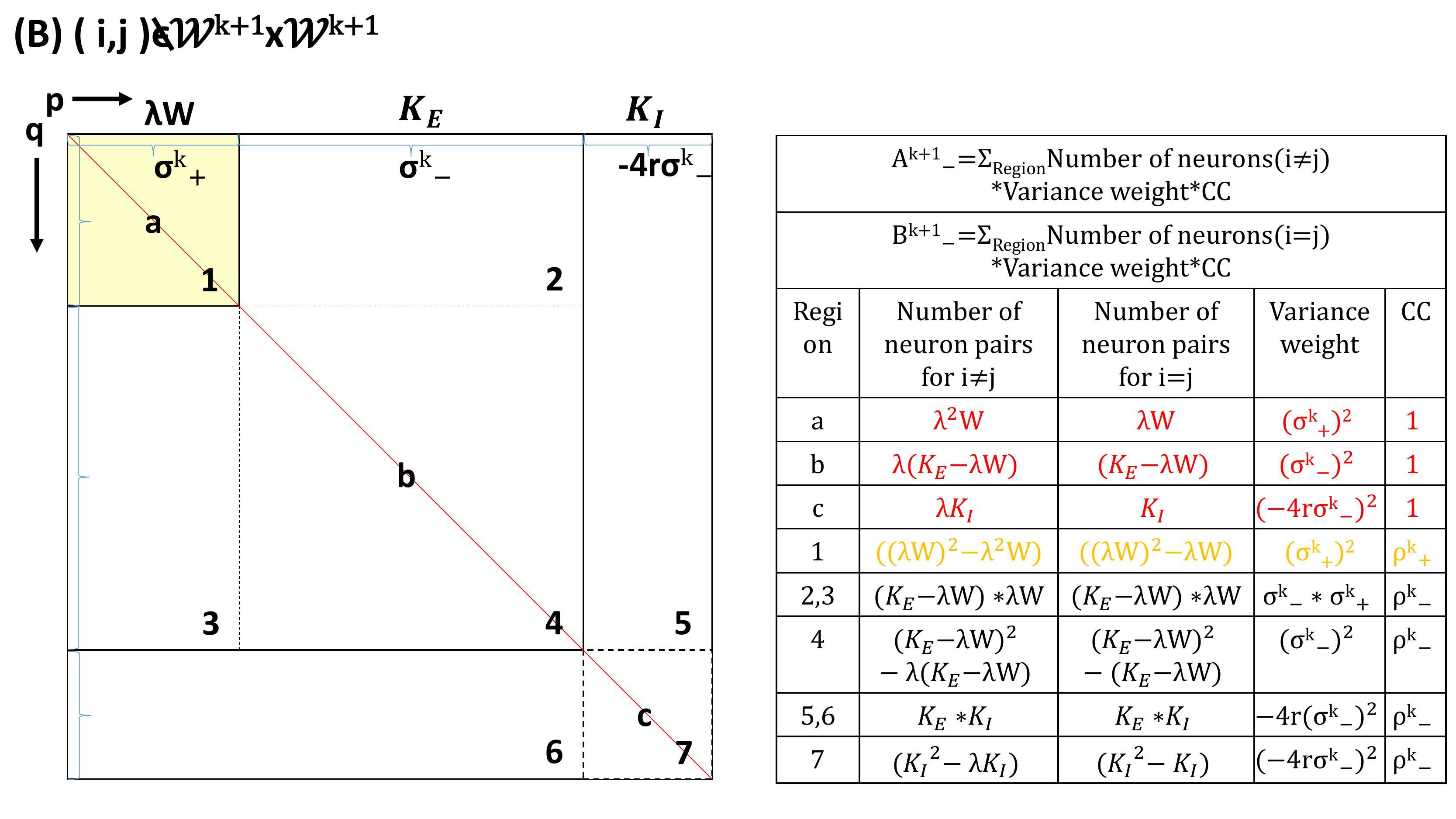}
\caption{The graphical computation of $\sum_{p,q}w_{ip}^{k+1}\sigma_{p}^{k}w_{jq}^{k+1}\sigma_{q}^{k}\rho^{k}_{pq}$ in  (\ref{Rho}). We lined up the subscripts $p=1,2,\dots,W,\dots,K_E,\dots,K_E+K_I$ and $q=1,2,\dots,W,\dots,K_E,\dots,K_E+K_I$ in two mutually perpendicular directions to form two squares for (A) $(i,j)\in\mathcal{W}^{k+1}\times\mathcal{W}^{k+1}$ and (B) $(i,j)\not\in\mathcal{W}^{k+1}\times\mathcal{W}^{k+1}$. Both of the squares can be decomposed into $3$ diagonal intervals and $7$ blocks according to the value of $\rho^{k}_{pq}$, the location of $p$ and $q$(in the packet or not) and the attribute of $p$ and $q$(excitatory or inhibitory). The value of number of neuron pairs, variance weight, correlation coefficient (CC) on each region and the expression of $A_{\pm}^{k+1}$ and $B_{\pm}^{k+1}$ are listed in the right tables.}
\label{fig2}
\end{figure*}
\section{Results}

\subsection{Non-constant coefficient of variation}
\begin{figure*}[t]
\begin{center}
\includegraphics[width=0.85\textwidth]{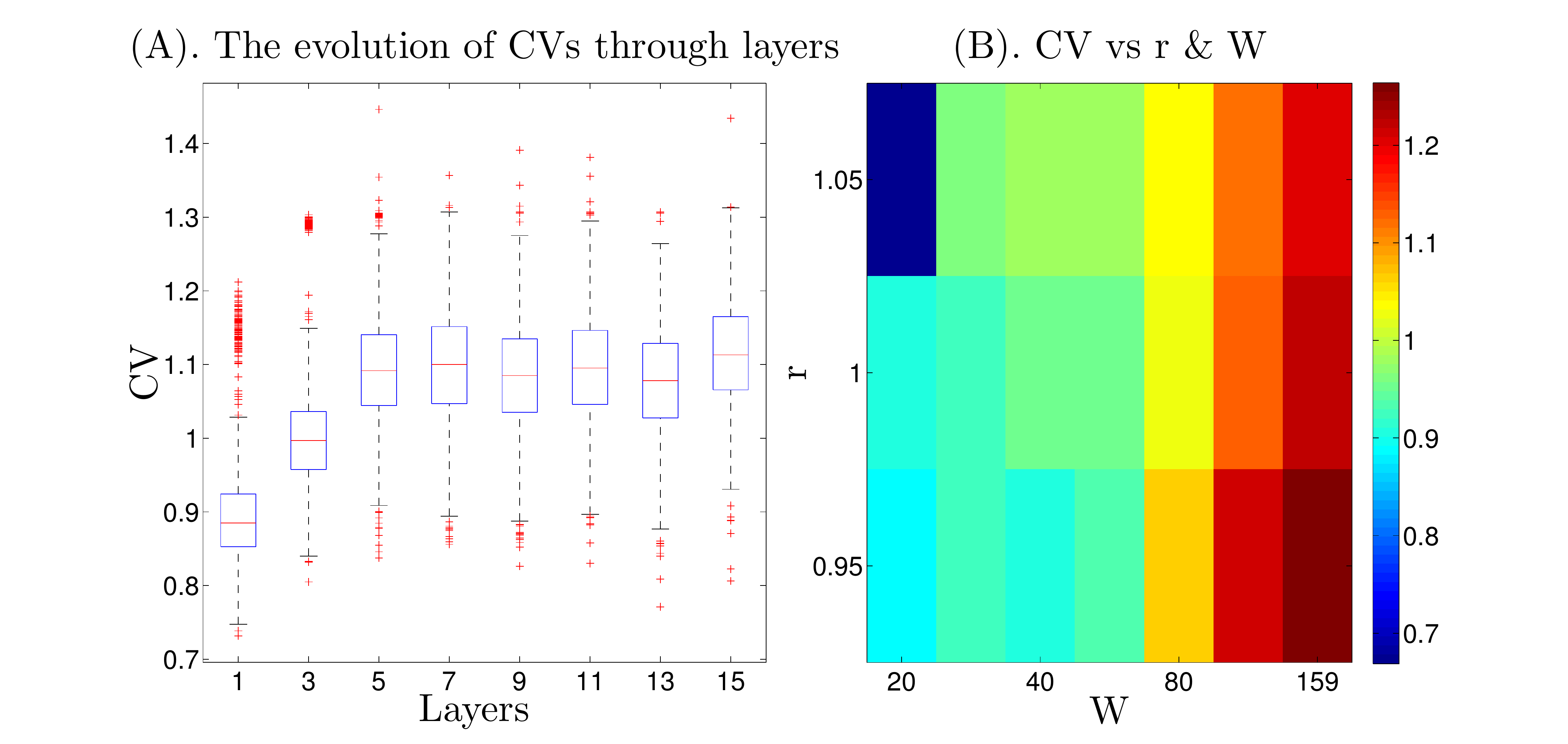}% Here is how to import EPS art
\caption{(A). Evolution of CVs of $2000$ LIF neurons in a FNN with respect to layer. (B). Averaged CV over 2000 neurons varies with parameters: synfire chain size $W$ and I-ESPS proportion $r$. 
}\label{fig3}
\end{center}
\end{figure*}

To our best knowledge, in most existing literature, the spike trains of LIF neurons were simulated by Poission or sub-/supra- Poisson processes through FNNs, and the coefficient of variation (CV) of the spike trains were assumed constant \cite{Aviel03,Rieke97}. Thus, the first order moment (mean) can be utilized to depict the second-order dynamics of synfire propagation.
However, this assumption could be away from the facts. To justify it,
we simulate $2000$ LIF neurons with initial Poisson input in synfire propagation to see the evolution of CV, which was defined in (\ref{CV}). As shown in Fig. \ref{fig3}(A), CV is non-stationary for layers. 
Moreover, even if the CVs go steady after many layers, they vary for parameters. For instance, there is a moderate rise of CVs with increasing cluster size $W$ or decreasing the proportion of ISPS and ESPS, denoted by $r$,  as shown in Fig. \ref{fig3}(B). This motivates us to employ Gaussian process to approximate the spike trains, instead of assuming fixed CV.
\subsection{Balanced network}
\begin{figure*}[!htb]
\begin{center}
\includegraphics[width=0.9\textwidth]{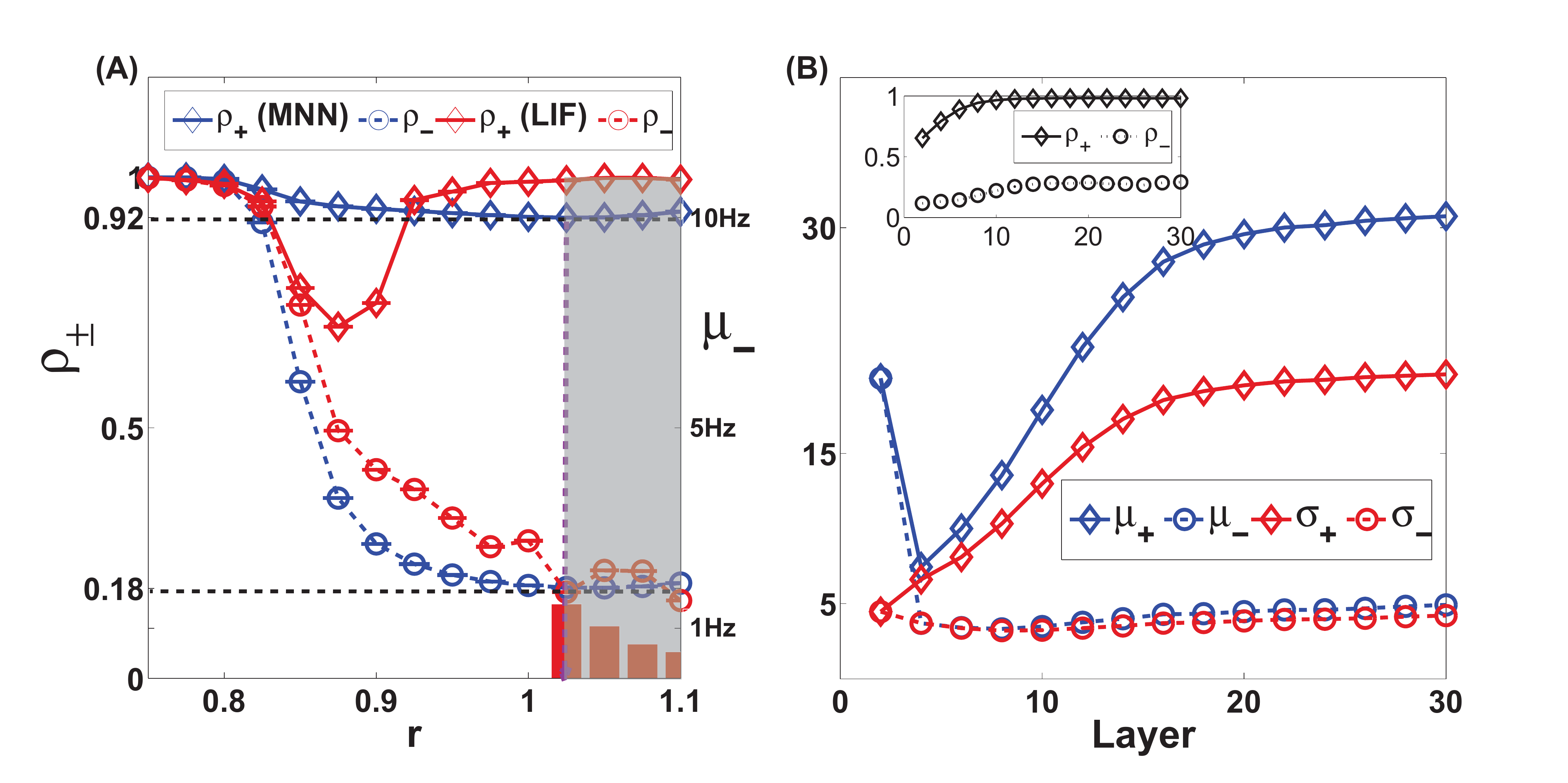}
\caption{Critical window of $r$ for the existence of synfire chain from the MNN model (mean, variance and correlation maps) and numerical LIF simulations. (A). Comparison between the theoretical results of the correlation map and simulation results of LIF neurons. $\rho_{+}$ (solid and $\diamond$) and $\rho_{-}$ (dash and $\circ$) vary with respect to the I-ESPS ratio $r$ with $W=32$ under the theoretical model (blue) or LIF simulation (red). $\rho_{\pm}$ are calculated at the $20$-th layer of the FNN. The mean firing rate (red bar and left vertical axis) of the out-of-packet neurons  {\bf $\mu_{-}$}, vary with respect to $r$ when $r>1$.
(B). Simulation results of the propagation of $\mu_{\pm}$, $\sigma_{\pm}$ and $\rho_{\pm}$ (the embedded subfigure) through layers by LIF network (\ref{LIF}). At the initial layer,  $N$ Possionian spiking trains are generated in a time interval of $20$ sec with the firing rate of
$20$Hz, among which $W$ trains in the packet are strongly correlated ($\rho\approx 1$) that represent and the rest out-of-packet $N-W$ trains are weakly correlated ($\rho< 0.2$).
}\label{fig4}
\end{center}
\end{figure*}

The parameter $r$ measures the ratio of ISPS over ESPS. $r=1$  means that the excitatory and inhibitory are exactly balanced.  Consider the correlation map in the mean field sense, (\ref{rho_in})-(\ref{rho_out}) and their expressions in detail. Noting $K_{E}=\lambda N_{E}$ with $N_{E}=4n=\frac{4}{5}N$, the terms of the highest order in the denominator, $B_{-}^{k+1}$, and the numerator, $A_{-}^{k+1}$, are both $(r-1)^2{K}_{E}^2(\sigma_{-}^k)^2\rho_{-}^k$ as $N\to\infty$. So, if $r\ne 1$, the identical term of the highest order in both denominator and numerator implies that $\rho_{-}^{k+1}=A_{-}^{k+1}/B_{-}^{k+1}$ approaches to $1$ as $k$ goes to infinity when $N$ is sufficiently large. Therefore, with a large size, $r=1$ is the critical value of the I-ESPS ratio for the existence of synfire chain in terms of property ${\mathbf P_{2}}$. As illustrated in Fig. \ref{fig4}(A), one can disclose that $\Lambda_{\epsilon,\delta}$ exists with some small values of $\epsilon$ and $\delta$ only in a balanced network, namely, $r\approx 1$, which implies the largest synfire gap ($\rho_{+}-\rho_{-}>0.7$). This finding can be verified by simulating LIF neuron in a feedforward network, as shown in Fig.\ref{fig4}(A) as well. However, when $r>1$, the out-of-packet mean firing rate  decreases quickly so that the out-of-packet spike frequency almost disappears (lower than $1~Hz$) when $r>1.2$. Therefore, to maintain a stable pattern of synchronization, the balanced network, namely, $r=1$, is necessary.
\subsection{Synfire chain condition and packet sizes}
\begin{figure*}[!htb]
\begin{center}
\includegraphics[width=0.8\textwidth]{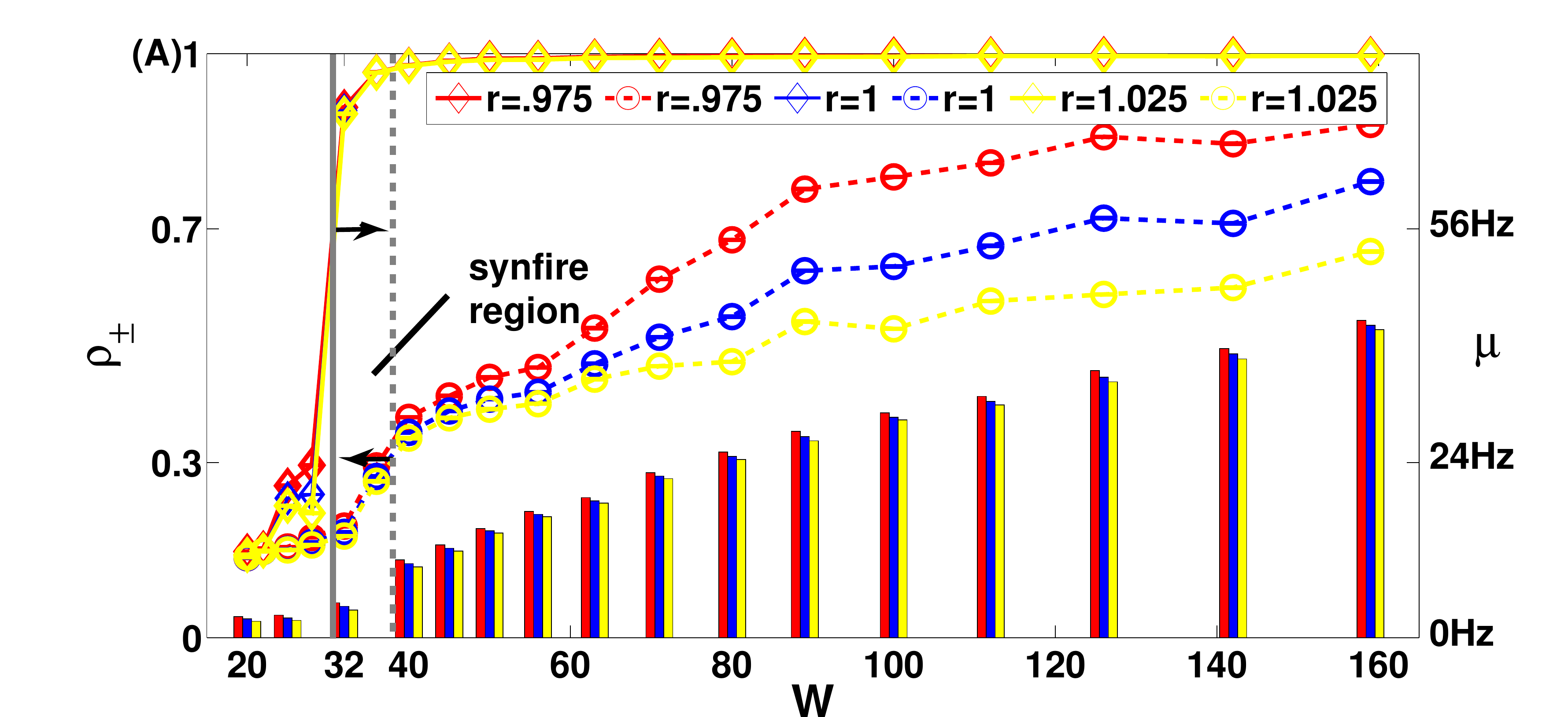}
\includegraphics[width=0.8\textwidth]{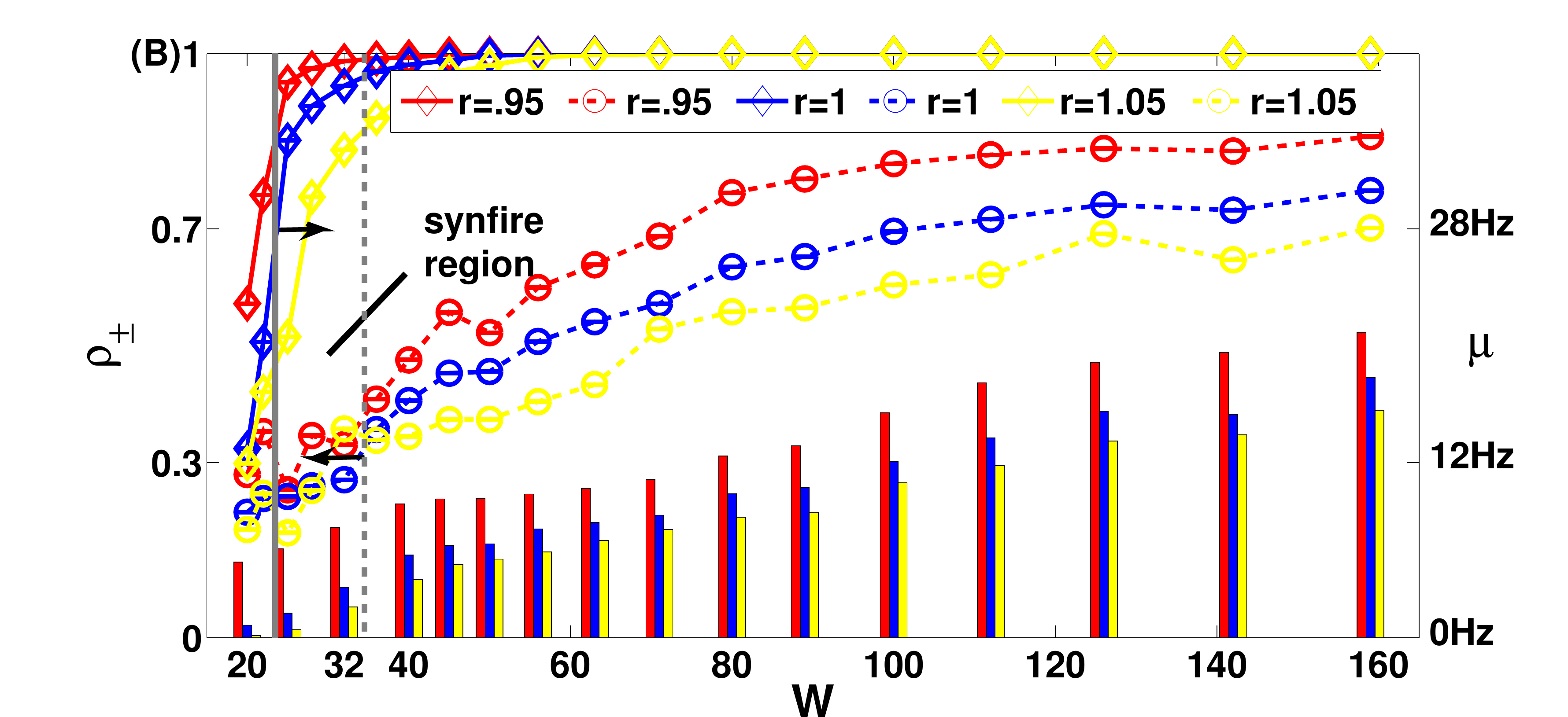}
\caption{The mean in-packet and out-packet correlation coefficients: $\rho_{+}$ (solid lines) and $\rho_{-}$ (dash lines) associated with the left vertical axis, and the mean firing rates $\mu$ (bars) associated with the right vertical axis vary with respect to packet size $W$, for different values of $r=0.95$ (red) for simulation (or $0.975$ for theoretical model), $1$ (blue) and $1.05$ (yellow) (or $1.025$ for theoretical model).  The synfire region (the shadow area) is defined by $\Lambda_{\epsilon,\delta}$ with $\epsilon=\delta=0.3$. The curves and bars are plotted for the theoretical model (A) and LIF network (B). The solid vertical line is for the critical annihilation packet size $W_a$ and the dotted line for invasion size $W_i$. (A). In the theoretical model, $W_a=31$ and $W_i=40$;  (B). In the LIF model, $W_{a}=24$ and $W_{i}=34$.}\label{fig5}
\end{center}
\end{figure*}

We discuss the packet size in the synfire chain:
$\Lambda_{\epsilon,\delta}$ with $\epsilon=\delta=0.3$ as an attracting set of the correlation map (\ref{Rho}) in a balanced network \cite{Brunel00,Vreeswijk98}. As illustrated in Fig. \ref{fig5}(A), to maintain $\Lambda_{\epsilon,\delta}$,  the packet of synfire chain should have an appropriate size. A large $W$ will enhance synchronization between neurons not only in-packet but also out-of-packet, and thus the in-packet synchrony will invade over the whole network, which leads the synfire gap to disappear due to the large packet size. That is, property ${\mathbf P}_{2}$ fails to hold (invasion). There exists a corresponding threshold value of $W$ by which synchronization invading occurs, named invasion packet size and denoted by $W_{i}$. In comparison, a very small $W$ may depress synchronization propagation in the packet and at the same time destruct the maintaining of the mean firing rates of the whole network, and thus the in-packet synfire chain will annihilate, which also leads the synfire gap disappears due to the small packet size. That is, property ${\mathbf P_{1}}$ cannot hold (annihilation). There also exists a corresponding threshold value of $W$, by which synchronization annihilation occurs, named annihilation  packet size and denoted by $W_{a}$. Therefore, an ideal $W$ to maintain a stable synfire chain is medium, belonging to the interval $(W_{a},W_{i})$, which leads both synchronous synfire packet and an asynchronous background with regular firing rates. As illustrated by Fig. \ref{fig5}(B), the simulation result by LIF neurons in the feedforward network has a good qualitative agreement with the theoretical results.

With the methods introduced in Section 2, we can derive an analytic result (a sufficient condition) with respect to $\epsilon$, $\delta$ and $\lambda$ for the existence of attractor $\Lambda_{\epsilon,\delta}$ in a balanced network ($r=1$). We consider the asymptotic stationary state of system (\ref{rho_in})-(\ref{rho_out}), in which all the mean-field variables such as $\sigma_{\pm}^t$, $\rho_{\pm}^t$, $A_{\pm}^t$ and $B_{\pm}^t$ have reached their steady states as shown in Fig. \ref{fig4}(B), whose values are denoted by $\sigma_{\pm}$, $\rho_{\pm}$, $A_{\pm}$ and $B_{\pm}$ respectively. For simplicity, let $P=5K_E\sigma_{-}^2(1-\rho_{-})$, $Q=[\sigma_{+}^2(\rho_{+}-\rho_{-})+(\sigma_{+}-\sigma_{-})^2\rho_{-}]$ and $R=2K_E\sigma_{-}(\sigma_{+}-\sigma_{-})\rho_{-}+(\sigma_{+})^2(1-\rho_{+})-(\sigma_{-})^2(1-\rho_{-})$.
With $r=1$, according to (\ref{rho_in})-(\ref{Variance_14}), by some algebras, we have
\begin{equation}
\rho_{+}\ge\frac{\lambda{P}+QW^2+RW}{P+QW^2+RW}
    \ge\frac{\lambda{P}+\lambda{Q}W^2+RW}{P+\lambda{Q}W^2+RW}
\end{equation}
and
\begin{eqnarray}
    |\rho_{-}|\le\frac{\lambda{P}+\lambda^2QW^2+\lambda{R}W}{P+\lambda^2QW^2+\lambda{R}W}.
\end{eqnarray}
Thus, it can be seen that if
\begin{eqnarray}
\frac{\lambda{P}+\lambda{Q}W^2+RW}{P+\lambda{Q}W^2+RW}>1-\epsilon,
\label{ineqXX}
\end{eqnarray}
and
\begin{eqnarray}
\frac{\lambda{P}+\lambda^2QW^2+\lambda{R}W}{P+\lambda^2QW^2+\lambda{R}W}<\delta,
\label{ineqXXX}
\end{eqnarray}
hold for some $W$, then the existence of $\Lambda_{\epsilon,\delta}$ can be guaranteed in the mean-field sense. By some algebras, one can derive that if
\begin{eqnarray}
\epsilon>\lambda\frac{1-\delta}{\delta}\label{syn_cond}
\end{eqnarray}
holds, named {\it synfire condition}, then (\ref{ineqXX})-(\ref{ineqXXX}) hold.

We are to identify the values of $W$ that satisfied (\ref{ineqXX})-(\ref{ineqXXX}) under the synfire condition. In fact, this synfire condition (\ref{syn_cond}) is equivalent to $\frac{1-\lambda}{\epsilon}-1<\frac{\delta-\lambda}{\lambda(1-\delta)}$. Thus, for any $\tau\in(\frac{1-\lambda}{\epsilon}-1,\frac{\delta-\lambda}{\lambda(1-\delta)})$, picking $W$ as the solution of equation $\frac{\lambda^2Q S^2+\lambda{R}S}{\lambda{P}}=\tau$ with respect to $S$, one can easily verify that (\ref{ineqXX})-(\ref{ineqXXX}) hold. Hence, we can derive a region of synfire size:
\begin{equation}
\begin{array}{r}
\operatorname{Reg}=\left\{W\right. \text { with } \frac{\lambda^{2} Q W^{2}+\lambda R W}{\lambda P}=\tau: \\
\left.\tau \in\left(\frac{1-\lambda}{\epsilon}-1, \frac{\delta-\lambda}{\lambda(1-\delta)}\right)\right\}.
\end{array}\label{reg}
\end{equation}
For each $W\in \operatorname{Reg}$, $\Lambda_{\epsilon,\delta}$ is an attracting set of the model (\ref{cc}) in the mean-field sense. In addition, $\max \operatorname{Reg}$ and $\min \operatorname{Reg}$ give the lower bound of the invasion packet size and the upper bound of the annihilation packet size respectively.
\subsection{Sparse synaptic density}
\begin{figure*}[!htb]
\begin{center}
\includegraphics[width=.6\textwidth]{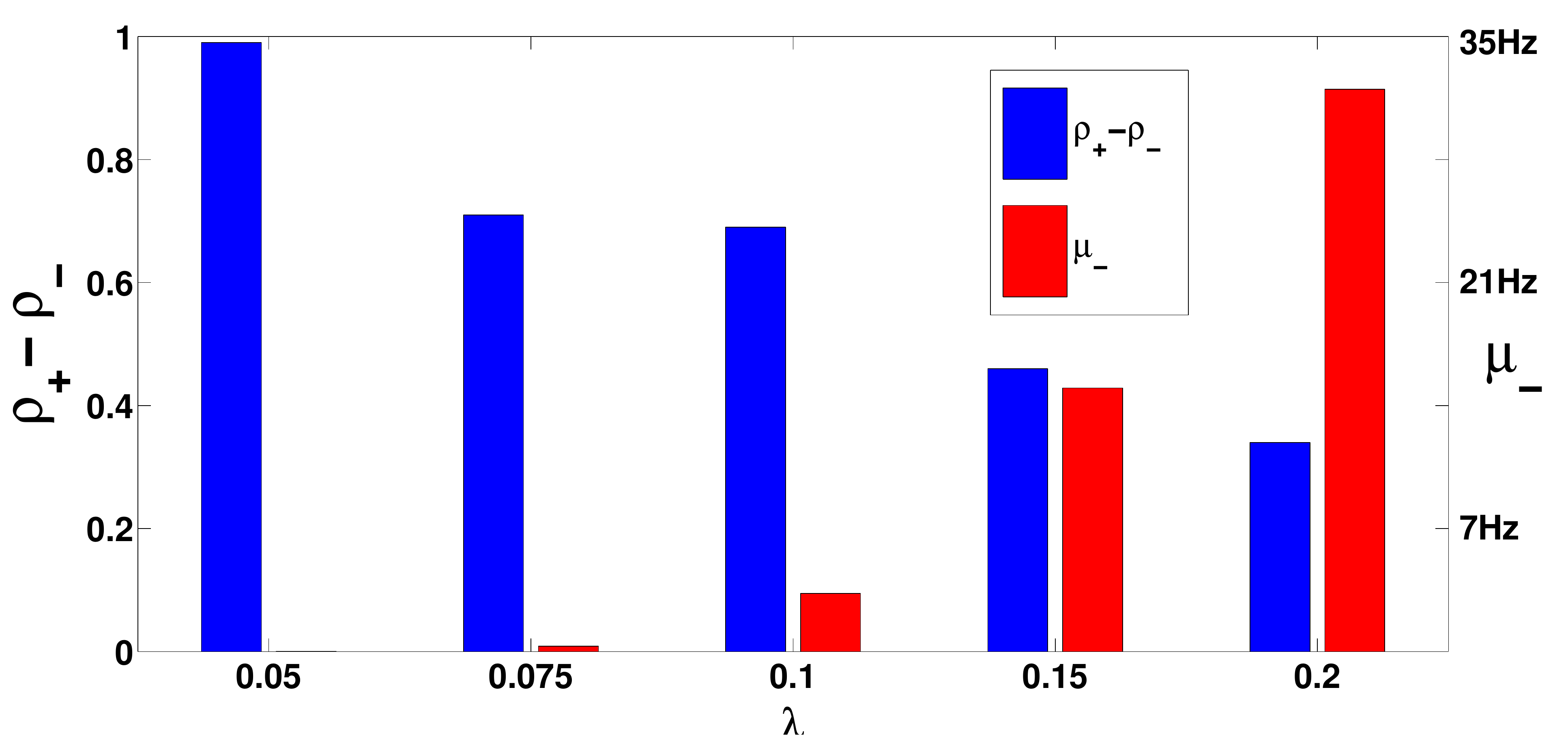}% Here is how to import EPS art
\caption{Barplots of the maximum synfire gap $\max_{W}\rho_{+}-\rho_{-}$ (blue) and the corresponding mean out-of-packet firing rate $\mu_{-}$ (red) with respect to linking density $\lambda$ in balanced LIF network.}\label{fig6}
\end{center}
\end{figure*}

Furthermore, the synfire condition (\ref{syn_cond}) reveals the dependence of synfire propagation on the synaptic density $\lambda$. The synfire condition (\ref{syn_cond}) and the region of synfire chain packet (\ref{reg}) imply that a smaller $\lambda$ leads a larger interval of synfire region in terms of the existence of $\Lambda_{\epsilon,\delta}$, which may cause a larger value of the maximum synfire gap $\max_{W}(\rho_{+}-\rho_{-})$. However, since the mean map in (\ref{MNN}) can be written as $\mu_{\pm}^{k+1}=\mathcal{S}_{1}(\hat{\mu}_{\pm}^{k+1},\hat{\sigma}_{\pm}^{k+1})$, where $\mathcal{S}_{1}(\cdot,\hat{\sigma}_{\pm}^{k+1})$ is a sigmoid function given $\hat{\sigma}_{\pm}^{k+1}$, and $\hat{\mu}_{-}^{k+1}=\lambda{W}(\mu_{+}^{k}-\mu_{-}^{k})+K_{E}(1-r)\mu_{-}^{k}$, in a balanced network ($r=1$), the mean out-of-packet firing rate will converge to $0$ through layers, because $\hat{\mu}_{-}^{k+1}\rightarrow0$ as $\lambda\rightarrow0$. That is to say, very small $\lambda$ results in the mean out-of-packet firing rate going extremely low, even near zero.
As shown in Fig. \ref{fig6}, the maximum synfire gap $\max_{W}\rho_{+}-\rho_{-}$ decreases with $\lambda$ but
an extremely small $\lambda$ (less than $0.075$) makes the spiking activities out of the packet vanish ($<0.3$ Hz). Therefore, to sustain a synfire chain and the asynchronous out-of-packet spiking activity , $\lambda$ should be taken a ``modestly'' small value. For instance, $\lambda=0.1$ was suggested in Refs.~\cite{Kumar10,Vogels05}.

\subsection{Synfire chain in recurrent neural network}
\begin{figure*}[!htb]
	\begin{minipage}[t]{0.5\textwidth}
		\centering
		\includegraphics[scale=0.8]{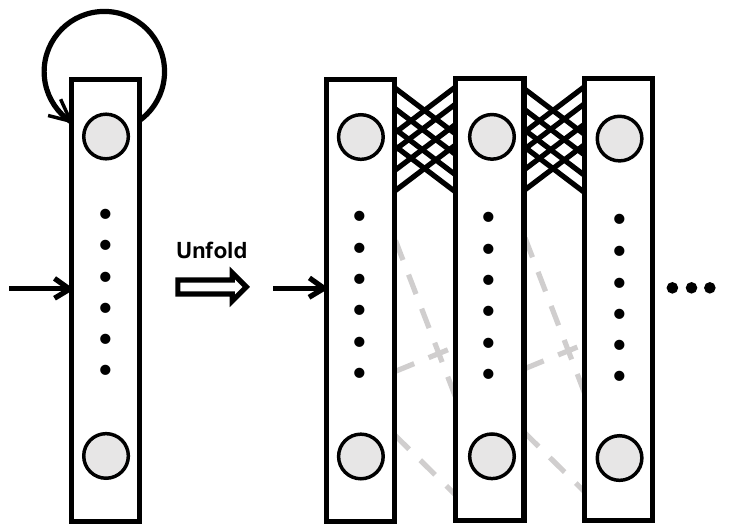}
		\caption*{(A)}
	\end{minipage}
	\begin{minipage}[t]{0.5\textwidth}
		\centering
		\includegraphics[scale=0.4]{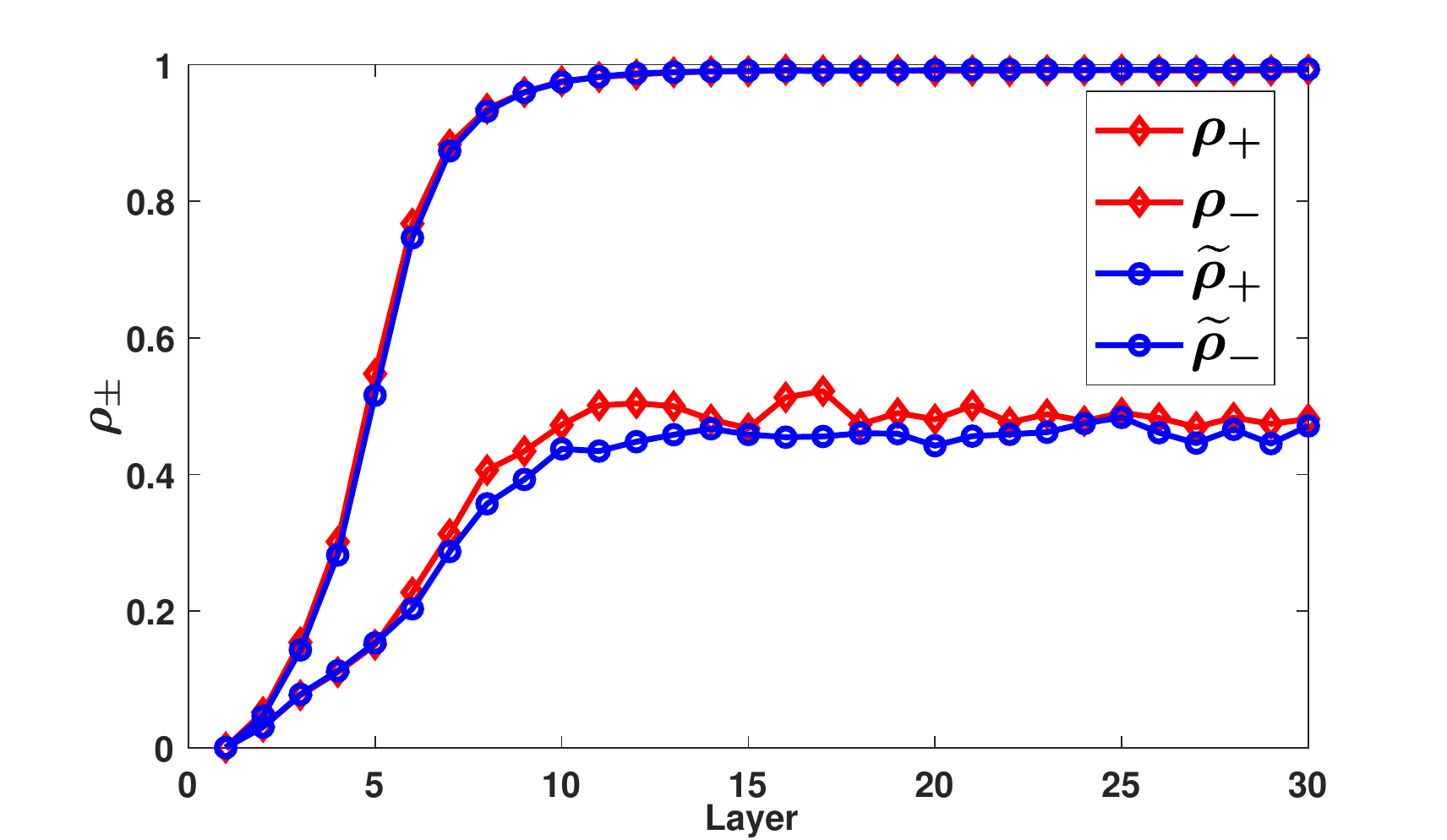}
		\caption*{(B)}
	\end{minipage}
	\caption{(A) The structure of recurrent spiking neural network and the equivalent form unfolded in discrete time. (B) The mean in-out packet correlation of recurrent ($\widetilde{\rho}_{\pm}$) and forward ($\rho_{\pm}$) LIF networks  with the same Poission input (20Hz, 20 seconds) and other parameters ($W=50$, $r=1$) in ten runs.}\label{fig7}
\end{figure*}

Besides FNN, we claim that above analysis still works in recurrent neural network (RNN)\cite{Lu2010}. If we consider the dynamics of a single layer of neurons in discrete time (see Fig. \ref{fig7}(A)), then (\ref{LIF}) can be written as:
\begin{eqnarray}\label{LIF2}
\tau_{m}dV_{i}^{s}(t)=- V_{i}^{s}(t) dt + I_{ext,i}^{s}+I_{syn,i}^{s},
\end{eqnarray}
where $s$ denotes time index. Similar to the Section 2.2, we iterate moment maps and get the moment closure:
\begin{eqnarray}
\begin{cases}
\mu^{s}(x)=\mathcal S_{1}^{x}\left[\mu^{s-1}(\cdot),cov^{s-1}(\cdot,\cdot)\right],\\
\sigma^{s}(x)=\mathcal S_{2}^{x}\left[\mu^{s-1}(\cdot),cov^{s-1}(\cdot,\cdot)\right],\\
\rho^{s}(x,y)=\Psi^{x,y}\left[\mu^{s-1}(\cdot),cov^{s-1}(\cdot,\cdot)\right],
\end{cases}\label{RMNN}
\end{eqnarray}
where the mapping functions are still defined as (\ref{momentmap}) and (\ref{cc}). Therefore FNN and RNN share almost the same evolution equations in our framework. For simplicity to compare, we denote time index $t$ in RNN as $k$ in the following paragraphs.

As shown in Fig. \ref{fig7}(A), to analyze RNN, firstly, we should unfold the network structure in discrete time. Then, RNN can be treated as FNN with constrained weight\cite{Mozer89}, which means packets $\mathcal W^k~(k=1,\cdots)$ share the same position and other random sparse connects are also fixed in different layers (weight sharing). If other parameters (e.g., I-ESPS proportion $r$, synaptic density $\lambda$, packet size $W$) take the same value as above, the dynamics of RNN and FNN are almost identical since we only focus on the forward propagation and do not consider the backpropagation like supervised learning. We denote the symbols using in FNN (Section 2.3) as $\widetilde{A}_{\pm}^{k}$, $\widetilde{B}_{\pm}^{k}$, $\widetilde{\rho}_{\pm}^{k}$ in RNN so that (\ref{expofAB}) turns to:
\begin{eqnarray}
\begin{cases}
\widetilde{A}_{+}^{k+1}=\left\langle\sum_{p,q}w_{ip}\widetilde{\sigma}_{p}^kw_{jq}\widetilde{\sigma}_{q}^k\widetilde{\rho}^{k}_{pq}\right\rangle_{i,j\in\mathcal W},\\
\widetilde{B}_{+}^{k+1}=\left\langle\sum_{p,q}w_{ip}\widetilde{\sigma}_{p}^kw_{iq}\widetilde{\sigma}_{q}^k\widetilde{\rho}^{k}_{pq}\right\rangle_{i\in\mathcal W},\\
\widetilde{A}_{-}^{k+1}=\left\langle\sum_{p,q}w_{ip}\widetilde{\sigma}_{p}^kw_{jq}\widetilde{\sigma}_{q}^k\widetilde{\rho}^{k}_{pq}\right\rangle_{i,j\notin\mathcal W},\\
\widetilde{B}_{-}^{k+1}=\left\langle\sum_{p,q}w_{ip}\widetilde{\sigma}_{p}^kw_{iq}\widetilde{\sigma}_{q}^k\widetilde{\rho}^{k}_{pq}\right\rangle_{i\notin\mathcal W}.\label{RexpofAB}
\end{cases}
\end{eqnarray}
Assuming that the input spikes in RNN and FNN are identical, with the mean field approximation, it is easy to get:
\begin{eqnarray}
\begin{cases}
\widetilde{A}_{\pm}^{k}=A_{\pm}^{k},\\
\widetilde{B}_{\pm}^{k}=B_{\pm}^{k},\label{ABinRNN}
\end{cases}
\end{eqnarray}
which leads to $\widetilde{\rho}_{\pm}^{k}={\rho}_{\pm}^{k}$. From Fig. \ref{fig7}(B), we can find that both networks have almost the same attractors. Furthermore, other dynamical characteristics referred in Section 3.1-3.4. can also be varified in the same way.

\section{Conclusion}
We have developed a theoretical framework of Gaussian random field to study how synchronization pattern propagates in feedforward or recurrent spiking neural network. Combined with the mean field approach, we analytically proved that the balance network is necessary for the stability of synfire chain in terms of suppressing the outside correlation and maintaining the in-packet spiking activities. And we derived a sufficient condition for a  stable synfire chain by providing an estimation of the appropriate packet size region which avoids both invasion and annihilation of synchronization, and revealing the role of spatial synaptic structure to sustain this spiking spatio-temporal pattern. Our analytic results show good agreements with the simulations of LIF network. We highlight that this approach based on moment closure is powerful and general to investigate propagation and stability of spatio-temporal patterns in random point field. One step further, we have not included the complex interactions like brain manifold structure in the random field here, as developed in Ref.\cite{Adler}, which is our future work and such a theoretical framework should be a valuable tool for investigating the varied dynamics of the spike patterns observed in experiments.

\section*{Acknowledgments}
\noindent This work is jointly supported by the National Key R\&D Program of China (No. 2019YFA0709502), the National Natural Sciences Foundation of China under Grant (No. 62072111), the 111 Project (No. B18015), the Shanghai Municipal Science and Technology Major Project (No. 2018SHZDZX01), ZJ LAB and the Shanghai Center for Brain Science and Brain-Inspired Technology.

\section*{Appendix A. Derivation of (\ref{Variance_11})-(\ref{Variance_14})}
\appendix
\renewcommand{\theequation}{A.\arabic{equation}}
As illustrated in the left squares of Fig. \ref{fig2}, we decompose the whole square into 10 regions according to the value of $\rho^{k}_{pq}$, the location of $p$ and $q$ (in the packet or not) and the attribute of $p$ and $q$ (excitatory or inhibitory). To calculate $\sum_{p,q}w_{ip}^{k+1}\sigma_{p}^{k}w_{jq}^{k+1}\sigma_{q}^{k}\rho^{k}_{pq}$, let us take region a,b and c for instance provided with $(i,j)\in\mathcal{W}^{k+1}\times\mathcal{W}^{k+1}$. Since region a,b,c are all diagonal intervals, the correlation coefficient $\rho^{k}_{pq}$ equals $1$ if $p$ connects $q$. In region a, note that neurons are fully connected within the cluster, $i,j$ share the same $W$ neighbors from previous layers. That is to say, both $p$ and $q$ go through the same index set of length $W$. In region b, $p$ and $q$ go through different index set of length $K_E-W$, thus $p$ connects $q$ for $K_E-W$ times if $i=j$ otherwise $\lambda(K_E-W)$ times, in which $\lambda$ describes the connection probability out of the cluster. In region c, $p$ and $q$ go through different index set of length $K_I$, thus $p$ connects $q$ for $K_I$ times if $i=j$ otherwise $\lambda{K_I}$ times. As for the variance weight, in region a, both $p$ and $q$ are in the cluster and excitatory, thus the variance weight is $(\sigma^k_+)^2$; in region b, both $p$ and $q$ are out of the cluster and excitatory, thus the variance weight is $(\sigma^k_-)^2$; in region c, both $p$ and $q$ are out of the cluster and inhibitory, thus the variance weight is $(-4r\sigma^k_-)^2$. Therefore, the contribution of these region to $A_{+}^{k+1}$ and $B_{+}^{k+1}$ are $W(\rho^k_+)^2*1$, $\lambda(K_E-W)(\rho^k_-)^2*1$, $\lambda{K}_I(-4r\rho^k_-)^2*1$ and $W(\rho^k_+)^2*1$, $(K_E-W)(\rho^k_-)^2*1$, $K_I(-4r\rho^k_-)^2*1$ respectively. To sum up the components in each region, we have
\begin{small}
\begin{align}\label{Variance_21}
&A_{+}^{k+1}
\!\approx\!\bigg\{W(\sigma_{+}^k)^2\!+\!\lambda(K_{E}\!-\!W)(\sigma_{-}^k)^2\!+\!\lambda{K_{I}}(-4r\sigma_{-}^k)^2\bigg\}\!*\!1\!+\!\bigg\{(W^2\!-\!W)(\sigma_{+}^k)^2\bigg\}\!*\!\rho_{+}^k\nonumber\\
&+\bigg\{2(K_{E}-W)W\sigma_{-}^k\sigma_{+}^k+[(K_{E}-W)^2-\lambda(K_{E}-W)](\sigma_{-}^k)^2
+2K_{E}K_{I}(-4r)(\sigma_{-}^k)^{2}\nonumber\\
&+(K_{I}^2-\lambda{K_{I}})(-4r\sigma_{-}^k)^2\bigg\}*\rho_{-}^k\nonumber\\
&=W(\sigma_{+}^k)^2(1\!-\!\rho_{+}^k)\!+\!\lambda(K_{E}\!-\!W)(\sigma_{-}^k)^2(1\!-\!\rho_{-}^k)\!+\!4r^2
\lambda{K_{E}}(\sigma_{-}^k)^2(1\!-\!\rho_{-}^k)\!+\!\rho_{-}^k\bigg\{(r\!-\!1)^2\nonumber\\
&\ K_{E}^2(\sigma_{-}^k)^2\!+\!2K_{E}W\sigma_{-}^k(\sigma_{+}^k\!-\!\sigma_{-}^k)
\!+\!W^2(\sigma_{+}^k\!-\!\sigma_{-}^k)^2\bigg\}\!+\!W^2(\sigma_{+}^k)^2(\rho_{+}^k\!-\!\rho_{-}^k),
\end{align}
\begin{align}\label{Variance_22}
&B_{+}^{k+1}
\!\approx\!\bigg\{W(\sigma_{+}^k)^2\!+\!(K_{E}\!-\!W)(\sigma_{-}^k)^2\!+{\!K_{I}}(\!-\!4r\sigma_{-}^k)^2\bigg\}*1+\bigg\{(W^2-W)(\sigma_{+}^k)^2\bigg\}*\rho_{+}^k\nonumber\\
&+\bigg\{2(K_{E}-W)W\sigma_{-}^k\sigma_{+}^k+[(K_{E}-W)^2-(K_{E}-W)](\sigma_{-}^k)^2
+2K_{E}K_{I}(-4r)(\sigma_{-}^k)^2\nonumber\\
&+(K_{I}^2-{K_{I}})(-4r\sigma_{-}^k)^2\bigg\}*\rho_{-}^k\nonumber\\
&=W(\sigma_{+}^k)^2(1\!-\!\rho_{+}^k)\!+\!(K_{E}\!-\!W)(\sigma_{-}^k)^2(1\!-\!\rho_{-}^k)
\!+\!4r^2{K_{E}}(\sigma_{-}^k)^2(1\!-\!\rho_{-}^k)\!+\!\rho_{-}^k\bigg\{(r-1)^2\nonumber\\
&\ K_{E}^2(\sigma_{-}^k)^2
\!+\!2K_{E}W\sigma_{-}^k(\sigma_{+}^k\!-\!\sigma_{-}^k)\!+\!W^2(\sigma_{+}^k\!-\!\sigma_{-}^k)^2\bigg\}\!+\!W^2(\sigma_{+}^k)^2(\rho_{+}^k\!-\!\rho_{-}^k),\end{align}
\begin{align}\label{Variance_23}
&A_{-}^{k+1}
\!\approx\!\bigg\{\lambda^2W(\sigma_{+}^k)^2\!+\!\lambda(K_{E}\!-\!\lambda{W})(\sigma_{-}^k)^2\!+\!\lambda{K_{I}}(-4r\sigma_{-}^k)^2\bigg\}*1
+\bigg\{((\lambda{W})^2-\lambda^2W)\nonumber\\
&(\sigma_{+}^k)^2\bigg\}\!*\!\rho_{+}^k\!+\!\bigg\{2(K_{E}\!-\!\lambda{W})\lambda{W}\sigma_{-}^k\sigma_{+}^k\!+\![(K_{E}\!-\!\lambda{W})^2\!-\!\lambda(K_{E}\!-\!\lambda{W})](\sigma_{-}^k)^2
\!+\!2K_{E}K_{I}\nonumber\\
&(-4r)(\sigma_{-}^k)^2+(K_{I}^2-\lambda{K_{I}})(-4r\sigma_{-}^k)^2\bigg\}\rho_{-}^k\nonumber\\
&=\lambda^2W(\sigma_{+}^k)^2(1-\rho_{+}^k)\!+\!\lambda(K_{E}-\lambda{W})(\sigma_{-}^k)^2(1-\rho_{-}^k)
+4r^2\lambda{K}_{E}(\sigma_{-}^k)^2(1-\rho_{-}^k)+\rho_{-}^k\nonumber\\
&\!\bigg\{\!(r\!-\!1)^2K_{E}^2(\sigma_{-}^k)^2
\!+\!2\lambda{K}_{E}W\sigma_{-}^k(\sigma_{+}^k\!-\!\sigma_{-}^k)\!+\!\lambda^2W^2(\sigma_{+}^k\!-\!\sigma_{-}^k)^2\!\bigg\}\!\!+\!\lambda^2W^2(\sigma_{+}^k)^2(\rho_{+}^k\!-\!\rho_{-}^k),
\end{align}
\begin{align}\label{Variance_24}
&B_{-}^{k+1}
\approx\bigg\{\lambda{W}(\sigma_{+}^k)^2+(K_{E}-\lambda{W})(\sigma_{-}^k)^2+K_{I}(-4r\sigma_{-}^k)^2\bigg\}*1
+\bigg\{((\lambda{W})^2-\lambda{W})\nonumber\\
&(\sigma_{+}^k)^2\bigg\}\!*\!\rho_{+}^k\!+\!\bigg\{2(K_{E}\!-\!\lambda{W})\lambda{W}\sigma_{-}^k\sigma_{+}^k\!+\![(K_{E}\!-\!\lambda{W})^2\!-\!(K_{E}\!-\!\lambda{W})](\sigma_{-}^k)^2
\!+\!2K_{E}K_{I}\nonumber\\
&(-4r)(\sigma_{-}^k)^2+(K_{I}^2-K_{I})(-4r\sigma_{-}^k)^2\bigg\}\rho_{-}^k\nonumber\\
&\!=\!\lambda{W}(\sigma_{+}^k)^2(1\!-\!\rho_{+}^k)\!+\!(K_{E}\!-\!\lambda{W})(\sigma_{-}^k)^2(1\!-\!\rho_{-}^k)
\!+\!4r^2{K_{E}}(\sigma_{-}^k)^2(1\!-\!\rho_{-}^k)\!+\!\rho_{-}^k\!\bigg\{\!(r\!-\!1)^2\nonumber\\
&\lambda^2{N}_{E}^2(\sigma_{-}^k)^2
\!+\!2\lambda{K}_{E}W\sigma_{-}^k(\sigma_{+}^k\!-\!\sigma_{-}^k)
+\lambda^2{W}^2(\sigma_{+}^k-\sigma_{-}^k)^2\bigg\}+\lambda^2{W}^2(\sigma_{+}^k)^2(\rho_{+}^k-\rho_{-}^k).
\end{align}
\end{small}
\bibliography{reference}

\end{document}